\documentclass[10pt]{amsart}
\usepackage{latexsym}
\usepackage{amsmath}
\usepackage{amssymb}
\usepackage[all]{xy}

\numberwithin{equation}{section}
\newtheorem{thm}{Theorem }[section]
\newtheorem{prop}[thm]{Proposition}
\newtheorem{lem}[thm]{Lemma}
\newtheorem{claim}[thm]{Claim}
\newtheorem{defi}[thm]{Definition}

\newtheorem{exmp}[thm]{Example}
\newtheorem{rem}[thm]{Remark}
\newtheorem{que}[thm]{Question}
\newcommand{\C}{{\mathbb C}}
\newcommand{\R}{{\mathbb R}}
\newcommand{\Z}{{\mathbb Z}}
\newcommand{\Q}{{\mathbb Q}}

\begin{document}
\title{ A Hasse diagram for  rational toral ranks}
\author{Toshihiro YAMAGUCHI}
\date{}
\address{Faculty of Education, Kochi University, 2-5-1,Kochi,780-8520, JAPAN}
\email{tyamag@kochi-u.ac.jp}
\maketitle

\begin{abstract}
Let $X$ be  a simply connected 
CW complex  with finite rational cohomology.
For the finite quotient  set  of rationalized
orbit spaces of $X$  obtained by almost free toral actions, ${\mathcal T}_0(X)
=\{[Y_i] \}$, induced by an equivalence relation   based on rational toral ranks,
we  order as $[Y_i]<[Y_j]$ 
if there is a rationalized Borel fibration $Y_i\to Y_j\to BT^n_{\Q}$
for some $n>0$. 
It presents
a variation of almost free toral actions on $X$.
We consider about the   Hasse diagram ${\mathcal H}(X)$ of the poset  
${\mathcal T}_0(X)$,
which makes a based graph  $G{\mathcal H}(X)$,
with some examples.
Finally  we will try to regard  $G{\mathcal H}(X)$ as the 1-skeleton of a finite CW complex  ${\mathcal T}(X)$
with base point $X_{\Q}$.
\end{abstract}

\section{Introduction} 
Let $ r_0(X)$ be the {\it rational  toral rank}
 of a simply connected CW complex  $X$
of $\dim H^*(X;\Q )<\infty $,
 i.e., the largest integer $r$ such that an $r$-torus
 $T^r=S^1 \times\dots\times S^1$($r$-factors)  can act continuously
 on a CW-complex $Y$ in  the rational homotopy type of $X$
 with all its isotropy subgroups finite (almost free action) \cite{AP},  \cite{FOT}, \cite{H}.
Recall that the  rationalized  
Borel space of almost free toral action $(ET^n\times_{T^n}^{\mu}Y)_{\Q}$
is homotopy equivalent to the rationalization of the orbit space of $Y$ obtained by 
the action $\mu$.
 In  a work of V.Puppe (for example see \cite{P}),
we can see a Hasse diagram of 
 the cohomology algebras of the  fixed point sets of circle  actions on $X$,
which are  correspond to 
the  rationalized  
Borel spaces,
from a point of view of a deformation.
We are interested in a rational  variation of  
Borel spaces of toral  actions with  no-fixed point
and the aim of this note is giving a framework 
for such an approach
based on rational toral ranks.

Due to the rational homotopy theory of D.Sullivan,
rationalized fibrations are  equivalent to   Koszul-Sullivan(KS)  extensions. 
Remark that,
when we give certain  KS-extensions of the Sullivan minimal model $M(X)$ 
of $X$ \cite{FHT},
  the existences of  the  free toral  actions  on 
 complexes  in   the rational homotopy type of $X$,
whose Borel fibrations induce the KS-extensions,  are  rationally 
guaranteed  by a result of S.Halperin \cite[Proposition 4.2]{H} (see Proposition 2.1 below).
We denote 
the homotopy set of  rationalized  Borel spaces of almost free toral actions $\mu$
 on 
complexes $X_{\mu}$ in the rational homotopy type of $X$,
which are given by certain KS-extensions (see \S 2), 
by  ${\mathcal X}=\coprod_{n=0}^{r_0(X)}{\mathcal X}_n$ where
${\mathcal X}_n:=\{(ET^n\times_{T^n}^{\mu}X_{\mu})_{\Q}\}$
for $n>0$ and ${\mathcal X}_0:=\{ X_{\Q}\}$.
For two
elements  $Y_1:=(ET^m\times_{T^m}^{\mu_1}X_1)_{\Q}$ and 
$Y_2:=(ET^n\times_{T^n}^{\mu_2}X_2)_{\Q}$ of ${\mathcal X}$
for $m<n$,
we denote  $Y_1<Y_2$
if 
there is a rationalized Borel fibration $Y_1\to Y_2\to BT^{n-m}_{\Q}$,
which 
satisfies  the homotopy commutative diagram 

$$\xymatrix{X_{\Q}\ar[d]\ar@{=}[r]& X_{\Q}\ar[d]&\\
Y_1\ar[r]\ar[d]&Y_2\ar[d]\ar[r]&BT^{n-m}_{\Q}\ar@{=}[d]\\
BT^{m}_{\Q}\ar[r]^{{Bi_n^m}_{\Q}}&
BT^{n}_{\Q}\ar[r]& BT^{n-m}_{\Q}.
}$$
Here  we put $Y_1:=X_{\Q}$ if $m=0$.
Then  $({\mathcal X},<)$
is a strict partially ordered  set (poset).


\begin{defi}
We give an equivalence relation of  ${\mathcal X}$  by
$Y_1\sim Y_2$ when  $Y_1,Y_2\in {\mathcal X}_n$ for some $n$
and $r_0(Y_1)=r_0(Y_2)$. 
For the quotient set ${\mathcal T}_0(X):={\mathcal X}/_{\sim}=\{ P_i\}_i$, 
  we put   $P_i<P_j$
 if   there are elements $Y_i,Y_j\in {\mathcal X}$ such that   $[Y_i]=P_i$,
$[Y_j]=P_j$
and $Y_i<Y_j$
or if  there is an element $P_k\in {\mathcal T}_0(X)$ with   $P_i<P_k$ and $P_k<P_j$.
\end{defi}
 
Notice that even if $T^m$ acts almost freely  on $X$ and $r_0(X)=n(>m)$,
then there does not always exist an almost free action of $T^{n-m}$
on a  complex  in  the rational homotopy type of the Borel space  $ET^m\times_{T^m}X$.
For example, when $X=S^3\times S^3\times S^7$,
we obtain  $r_0(X)=3$ 
by standard $T^3$-action
 $(s_1,s_2,s_3)\cdot (z_1,z_2,z_3)=(s_1z_1,s_2z_2,s_3z_3)$.
But  there exists  a free $S^1$-action $\mu:S^1\times Y\to Y$ for a finite complex
$Y$ with $Y_{\Q}\simeq X_{\Q}$ and 
$r_0(ES^1\times_{S^1}^{\mu}Y)=0$.
It is also rationally given as the total space of a non-trivial fibration
with fiber 
$\C P^3$ and base $S^3\times S^3$.
See Example 3.5 below  for detail.
Thus  we stand on our starting point.

\begin{claim}
The poset ${\mathcal T}_0(X)=(\{P_i\}_i,<)$ is not  
totally ordered in general.
\end{claim}

The poset 
${\mathcal T}_0(X)$
makes a Hasse diagram of the   sets  $\{P_i\}_i$.
We denote it as  ${\mathcal H}(X)$.
It is not a numerical but is  a {graphical  (rational)  homotopy invariant} of  spaces.
Here we can put 
 $i<j$ if $P_i<P_j$ and  fix 
  $P_0=[X_{\Q}]$,
  $P_1=[(ES^1\times_{S^1}Y)_{\Q}]$ with $r_0(ES^1\times_{S^1}Y)=r_0(X)-1$,
$\cdots , $ $P_{r_0(X)}=[(ET^{r_0(X)}\times_{T^{r_0(X)}}Y)_{\Q}]$ 
 with $r_0(ET^{r_0(X)}\times_{T^{r_0(X)}}Y)=0$
for a complex $Y$ in the rational homotopy of $X$.
The  subset $\{P_1,..,P_{r_0(X)}\}$ always  exists  by the  restrictions to
$i^m_{r_0(X)}(T^m)=\{ (s_1,..,s_m,1,..,1)|s_i\in S^1\}$
of an almost free  
$T^{r_0(X)}$-action  on $Y$
for $m=1,..,r_0(X)$.
We observe  from the above definition

\begin{lem} 
(1) For $Y\in{\mathcal X}$, $r_0(Y)=n$ if and only if $n=max\{k\ |\ [Y]=P_i<P_{i_1}<\cdots <P_{i_k}, \ r_0(P_{i_k})=0\}$.
Then  the path 
of length $n$, $P_i\to P_{i_1}\to \cdots \to P_{i_n}$, is unique in the graph.
In particular, if $Y=X_{\Q}$,
 $P_i=P_0$, $P_{i_j}=P_j$ for $j=1,..,n$.\\
(2) For $Y\in{\mathcal X}$, $Y\in {\mathcal X}_n$ if and only if 
$n=max\{k\ |\ P_0<P_{i_1}<\cdots  <P_{i_k}=[Y]\}$.
Then $n=d(P_0, [Y])$, 
 the distance between $P_0$ and $[Y]$
in the  graph.
  \end{lem}

Our Hasse diagrams are restricted to certain forms.
 Of course,  ${\mathcal T}_0(X)$ is a  finite set.
Especially, when $r_0(X)=n$,
$n<\sharp {\mathcal T}_0(X)\leq (n+(n-1)+\cdots +2+1)+1=(n^2+n)/2+1$
 ($0<\sharp \ {\mathcal X}_m/_{\sim}\leq n-m+1$ for $m\leq n$).
In particular,   $\sharp {\mathcal T}_0(X)=r_0(X)+1$ if and only if 
  ${\mathcal T}_0(X)$ is totally ordered.
For example, if  $r_0(X)=3$ and $r_0(X')=4$
for some spaces $X$ and $X'$, $4\leq \sharp {\mathcal T}_0(X)\leq 7$ and 
$5\leq \sharp {\mathcal T}_0(X')\leq 11$
and 
 ${\mathcal H}(X)$ and ${\mathcal H}(X')$ are  certain 
sub-diagrams (see Remark 3.7 below) of  the Hasse diagrams:

$$
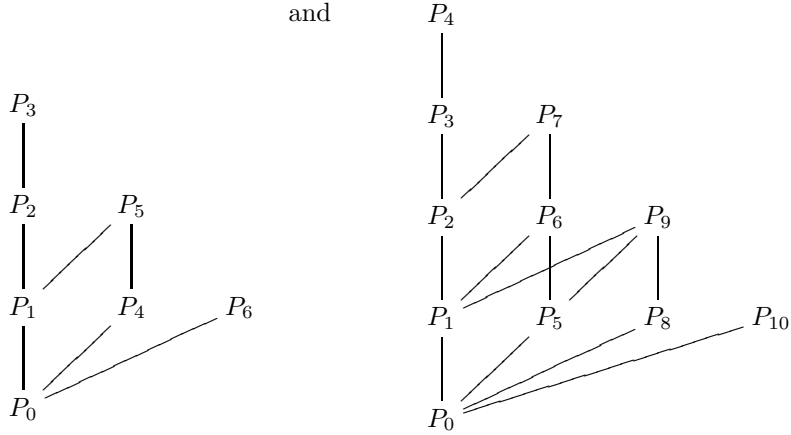
\xymatrix{ 
& & \\
P_3 &  &   \\
P_2 \ar@{-}[u]& P_5&  \\
P_1 \ar@{-}[u]\ar@{-}[ru] & P_4\ar@{-}[u]& P_6  \\
P_0\ar@{-}[u]\ar@{-}[ur]\ar@{-}[urr]\\
}\ \ \ \mbox{and} \ \ \ \ \ \ \ \ \  \ 
\xymatrix{ 
P_4 & &  & \\
P_3 \ar@{-}[u]& P_7 &  & \\
P_2 \ar@{-}[u]\ar@{-}[ru]& P_6 \ar@{-}[u]& P_{9}& \\
P_1 \ar@{-}[u]\ar@{-}[ru]\ar@{-}[urr] & P_5\ar@{-}[u]\ar@{-}[ru] & P_8\ar@{-}[u]& P_{10}  \\
P_0\ar@{-}[u]\ar@{-}[ur]\ar@{-}[urr]\ar@{-}[urrr]\\
}$$
, respectively.
If there exist such spaces, 
$r_0(P_0)=3$, 
 $r_0(P_1)=2$, $r_0(P_2)=r_0(P_4)=1$, 
 $r_0(P_3)=r_0(P_5)=r_0(P_6)=0$ in the left hand and
$r_0(P_0)=4$, 
 $r_0(P_1)=3$, $r_0(P_2)=r_0(P_5)=2$, 
 $r_0(P_3)=r_0(P_6)=r_0(P_8)=1$,
 $r_0(P_4)=r_0(P_7)=r_0(P_9)=r_0(P_{10})=0$ in the right hand.
 Here $r_0(P_i)$ means $r_0(Y)$ for some space $Y$ with $P_i=[Y]$.

We can describe    a  point $P_i=[Y]$ of  
${\mathcal T}_0(X)$  by  the double index (lattice point)
$$d.i.(P_i):=(s,t)\ \ \ \ \  ; \ \  s+t\leq r_0(X)$$  when  
$$Y\in {\mathcal X}_t\mbox{ \ \  and \ \ } r_0(Y)=r_0(X)-s-t$$
 by Definition 1.1. 
If $P_i\neq P_j$
in ${\mathcal T}_0(X)$, $d.i.(P_i)\neq d.i.(P_j)$.
 For example, in the above right diagram of $r_0(X)=4$, we see $d.i.(P_0)=(0,0)$, 
$d.i.(P_1)=(0,1)$, $d.i.(P_2)=(0,2)$,  $d.i.(P_3)=(0,3)$, $d.i.(P_4)=(0,4)$, $d.i.(P_5)=(1,1)$, $d.i.(P_6)=(1,2)$,  $d.i.(P_7)=(1,3)$, $d.i.(P_8)=(2,1)$, $d.i.(P_9)=(2,2)$ and $d.i.(P_{10})=(3,1)$.
In general, when $r_0(X)>1$,
if there is a  circle action on $X$ that represents $P$ with $d.i.(P)=(r_0(X)-1,1)$,
then it is a ``bad'' action in a meaning
since the orbit space permits no almost free circle action.

\begin{claim}
(1)  If $P_i<P_j$ for $d.i.(P_i)=(s,t)$ and  $d.i.(P_j)=(s',t')$, then $s\leq s'$ and $t<t'$.\\
(2) If there is a point $P_i$ with $d.i.(P_i)=(s,t)$, then 
there are points $\{ P_j\}$ with double indexes $(s,t+1),.., (s,r_0(X)-s)$, too.
\end{claim}

Notice that 
a Hasse diagram ${\mathcal H}$ can be seen as   a connected, finite,
non-directed, simple 
graph $G{\mathcal H}$ with base point corresponding to the minimal element
in general.
We  say a graph with a base point as  a  {\it based graph} in this paper.
Define 
$\phi :{\mathcal T}_0(X)\to \Z_{\geq 0}\times \Z_{\geq 0}$
by $\phi (P):=d.i.(P)$ and extend 
$\tilde{\phi} :G{\mathcal H}(X)\to \R_{\geq 0}\times \R_{\geq 0}$
by $\tilde{\phi} (P_iP_j)=d.i.(P_i)-d.i.(P_j)$, the line segment with extremal points  $d.i.(P_i)$ and 
$d.i.(P_j)$. Then there is a commutative diagram 
 
$$\xymatrix{ {\mathcal T}_0(X)
\ar[r]^{\phi\ \ \ }\ar[d]_{\cap}& \Z_{\geq 0}\times \Z_{\geq 0}\ar[d]^{\cap}\\
G{\mathcal H}(X)\ar[r]^{\tilde{\phi}\ \ \ }&  \R_{\geq 0}\times \R_{\geq 0}.\\
}
$$

Note that 
$\tilde{\phi}$ is injective, that is, $\tilde{\phi}$ gives the realization of 
${\mathcal H}(X)$ into $\R_{\geq 0}\times \R_{\geq 0}$
induced by the above double indexes.
We see ${\mathcal H}(X)= {\mathcal H}(Y)$
if and only if 
$\tilde{\phi}G{\mathcal H}(X)= \tilde{\phi}G{\mathcal H}(Y)$.
On the other hand,
we can reconstruct $\tilde{\phi}G{\mathcal H}(X)$
from $G{\mathcal H}(X)$ graphically (see \S 4).
Thus

\begin{thm}
 For some spaces $X$ and $Y$,  ${\mathcal H}(X)= {\mathcal H}(Y)$  if and only if   
 $G{\mathcal H}(X)$ and $G{\mathcal H}(Y)$ are  isomorpic
as  based graphs.
\end{thm}

There do not exist the following Hasse diagrams in our ones:
{\small 
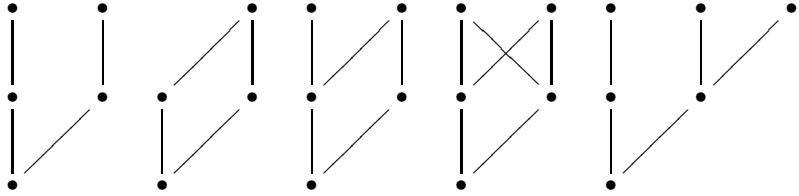
$$
\xymatrix{ 
\bullet &\bullet\\
\bullet \ar@{-}[u] & \bullet\ar@{-}[u]\\
\bullet\ar@{-}[u]\ar@{-}[ur]
}\ \ \ \ 
\xymatrix{ 
 &\bullet\\
\bullet \ar@{-}[ru] & \bullet\ar@{-}[u]\\
\bullet\ar@{-}[u]\ar@{-}[ur]
}\ \ \ \  
\xymatrix{ 
\bullet &\bullet\\
\bullet \ar@{-}[ru] \ar@{-}[u] & \bullet\ar@{-}[u] \\
\bullet\ar@{-}[u]\ar@{-}[ur]
}\ \ \ \  
\xymatrix{ 
\bullet &\bullet\\
\bullet \ar@{-}[ru] \ar@{-}[u] & \bullet\ar@{-}[u]\ar@{-}[lu] \\
\bullet\ar@{-}[u]\ar@{-}[ur]
}\ \ \ \  
\xymatrix{ 
\bullet &\bullet&\bullet\\
\bullet  \ar@{-}[u] & \bullet\ar@{-}[u]\ar@{-}[ru] &\\
\bullet\ar@{-}[u]\ar@{-}[ur]&
} $$}
\noindent
$\cdots$. 
 For example, the graph 
 $$\xymatrix{ 
A \ar@{-}[r] & B \ar@{-}[r] & C \ar@{-}[r] & D \ar@{-}[r] & E  \\
}$$
represents  the totally ordered Hasse diagram of a space with rational toral rank $4$ 
 if we choose the base point as  $A$ or $E$.
Also  it represents  (2)  of Example 3.5 if we choose the base point as  $B$ or $D$.
 But  if we choose the base point as  $C$,
 the graph  corresponds to non of our Hasse diagrams. 
  Also  the graph 
 $$\xymatrix{& F \ar@{-}[r]\ar@{-}[d]& G\ar@{-}[d] &&\\  
A \ar@{-}[r] & B \ar@{-}[r] & C \ar@{-}[r] & D \ar@{-}[r] & E  \\
}$$
represents  our   Hasse diagrams $(a)$ or $(b)$ below   
 if we choose the base point as  $A$ or $B$, respectively.
 
$$(a)\ \ 
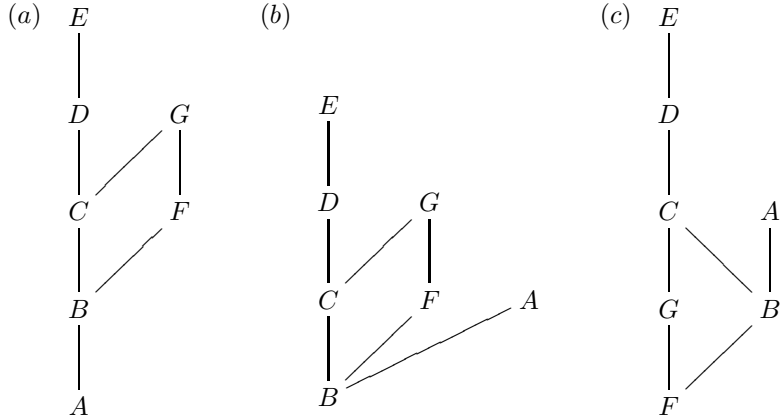
\xymatrix{ 
E &   \\
D \ar@{-}[u]& G \\
C \ar@{-}[u]\ar@{-}[ru] & F\ar@{-}[u]  \\
B\ar@{-}[u]\ar@{-}[ur]\\
A\ar@{-}[u]
}\ \ \ \ \  \ \ 
 (b) \ \ \xymatrix{ 
 & & \\
E &  &   \\
D \ar@{-}[u]& G&  \\
C \ar@{-}[u]\ar@{-}[ru] & F\ar@{-}[u]& A \\
B\ar@{-}[u]\ar@{-}[ur]\ar@{-}[urr]\\
}\ \ \ \  \ \ (c) \ \
 \xymatrix{ 
E&\\
D\ar@{-}[u] &    \\
C \ar@{-}[u]& A \\
G \ar@{-}[u] & B\ar@{-}[lu] \ar@{-}[u]\\
F\ar@{-}[u]\ar@{-}[ur]\\
}
$$
If we choose the base point as  $F$,
 the Hasse diagram is given as $(c)$, which is not ours
since the points $G$ and $B$ must be  a  same one
from Definition 1.1.
Also we can check that 
the other points are not impossible to be realized as the minimal elements of our Hasse diagrams
(the base point of  $G{\mathcal H}$). 
Note that the author does not know whether or not exists 
a space (rational model)  $X$ with ${\mathcal H}(X)=(a)$.
The following question is essential.

\begin{que}
Find an example of two  spaces (rational models)  $X$ and $Y$
such that $\phi {\mathcal T}_0(X)=\phi {\mathcal T}_0(Y)$ in $\Z_{\geq 0}\times \Z_{\geq 0}$
 but ${\mathcal H}(X)\neq {\mathcal H}(Y)$.
 \end{que}

\begin{rem}
Our  definition of $({\mathcal X}/\sim ,\ >)$ in  Definition 1.1 may be   rough.
But if we do not take the quotient,
the poset $({\mathcal X},>)$ seems  very complicated.
For example, even when $X=S^3\times S^3$,
the Hasse diagram is 

$$\xymatrix{ 
\cdots P_i \cdots   \\
(S^2\times S^3)_{\Q} \ar@{-}[u]_{\bf \cdots}^{\bf \cdots}\\
(S^3\times S^3)_{\Q} \ar@{-}[u]
}$$
(it seems as a broom)
where $H^*(S^2\times S^3;\Q)\cong \Q [t_1]/(t_1^2)\otimes \Lambda (v)$
with   $|v|=3$ and 
$$\{P_i\}_{i}\cong \{ H^*(P_i;\Q)\}_i\cong \{\frac{\Q [t_1,t_2]}{(t_1^2+a_it_2^2,t_1t_2)}
|\ a_i\in \Q^* \} \cong \Q^*/{\Q^*}^2,$$
which is an infinite set ($\Q^*=\Q-0$ is the unite group of $\Q$).
Note that 
$D_1u=t_1^2$,
$D_1v=0$,
$D_2u=t_1^2+a_it_2^2$ and 
 $D_2v=t_1t_2$  
for $M(X)=(\Lambda (u,v),0)$   (see \S 2). 
\end{rem}

\begin{rem}  { For $Y\in{\mathcal X}$, 
$ {\mathcal T}_0(Y)\equiv \{P_i\in {\mathcal T}_0(X) \ | \ [Y]=P_i\ or \ [Y]<P_i\}$ as  ordered sets.
Thus  ${\mathcal H}(Y)$
is a sub-Hasse diagram of ${\mathcal H}(X)$.
Also
for two spaces $X$ and $X'$, 
$G{\mathcal H}(X\times X')\supset  
G{\mathcal H}(X)\vee G{\mathcal H}(X')$
as a  subgraph   with  vertexes $\{P_i\}_i$ and edges $\{ P_iP_j\}=\{ P_i<P_j| \ \mbox{ there is 
 no } P_k\ \mbox{ with } P_i<P_k<P_j\}_{i,j}$.
Here the right hand is the one point union  $G{\mathcal H}(X)\coprod  G{\mathcal H}(X')/\sim$
where $P_{r_0(X)}\sim P_0'$ for ${\mathcal T}_0(X)=\{ P_i\}_i$
 and ${\mathcal T}_0(X')=\{ P'_i\}_i$. 
It is  a  grafting of one on the other.
By using a Sullivan model, 
 S.Halperin 
 indicates that rational  toral rank does not preserve the product fomula
$r_0(X\times X')=r_0(X)+r_0(X')$ in general  \cite{JL}(\cite[Ex.7.19]{FOT}).
Thus this embedding   may be  complicated in general   
(see Example 3.9 below).
}\end{rem}

\noindent
{\bf Acknowledgement}. 
The author would like to thank 
Katsuhiko Kuribayashi,  Shizuo Kaji and the referee
for their valuable suggestions and  is grateful  to Yves F\'{e}lix
for his encouragement.


\section{A Halperin's result}

Let $X$ be a simply connected CW complex  of finite type and  
the  Sullivan minimal model $M(X)=(\Lambda {V},d)$. 
  It is a free $\Q$-commutative differential graded algebra 
 with a $\Q$-graded vector space $V=\bigoplus_{i\geq 2}V^i$
 where $\dim V^i<\infty$ and a decomposable differential; i..e., $d(V^i) \subset (\Lambda^+{V} \cdot \Lambda^+{V})^{i+1}$ and $d \circ d=0$.
 Here  $\Lambda^+{V}$ is 
 the ideal of $\Lambda{V}$ generated by elements of positive degree. 
Denote the degree of a homogeneous element $x$ of a graded algebra as $|{x}|$.
Then  $xy=(-1)^{|{x}||{y}|}yx$ and $d(xy)=d(x)y+(-1)^{|{x}|}xd(y)$. 
Note that  $M(X)$ determines the rational homotopy type of $X$, $X_{\Q}$.
In particular,  $H^*(\Lambda {V},d)\cong H^*(X;\Q )$.
Refer  \cite{FHT} for detail.

If an $r$-torus $T^r$ acts on $X$
by $\mu :T^r\times X\to X$, there is the Borel fibration
$$
X \to ET^r \times_{T^r}^{\mu} X \to BT^r,
$$
where 
$ ET^r \times_{T^r}^{\mu} X $  is the orbit space of the  action
$g(e,x)=(eg^{-1},gx)$  
on the product $ ET^r \times  X $.
It is rationally given by the  KS extension (model)
$$
(\Q[t_1,\dots,t_r],0)
 \to (\Q[t_1,\dots,t_r] \otimes \Lambda {V},D)
 \to (\Lambda {V},d)\ \ \ \ (*)$$
where    with $|{t_i}|=2$ for $i=1,\dots,r$, $Dt_i=0$ and
$Dv \equiv dv$ modulo the ideal $(t_1,\dots,t_r)$ for $v\in V$.

\begin{prop}\cite[Proposition 4.2]{H}
Suppose that $X$ is a simply connected CW-complex  with 
$\dim H^*(X;\Q)<\infty$.
Put $M(X)=(\Lambda V,d)$.
Then  $r_0(X) \ge r$ if and only if there is a KS extension $(*)$
 satisfying $\dim H^*(\Q[t_1,\dots,t_r] \otimes \wedge{V},D)<\infty$.
 Moreover,
 if  $r_0(X) \ge r$,
 then $T^r$ acts freely on a finite complex $Y$
 that has  the same rational homotopy type as $X$
 and $M(ET^r\times_{T^r}Y)\cong 
(\Q[t_1,\dots,t_r] \otimes \wedge{V},D)$.
\end{prop}

Thus we can put
$${\mathcal X}_n=\{ (\Q[t_1,\dots,t_n] \otimes \Lambda {V},D) \ 
|\ \dim H^*(\Q[t_1,\dots,t_n] \otimes \wedge{V},D)<\infty\}/_{\cong} $$
 for $M(X)=(\Lambda V,d)$.
The KS extension of the fibration $Y_1\to Y_2\to BT^{n-m}_{\Q}$ in \S 1
is given by the homotopy commutative diagram 
$$\xymatrix{ & (\Lambda V,d)\ar@{=}[r]&(\Lambda V,d)\\
(\Q [t_{m+1},..,t_{n}],0)\ar[r]& (\Q [t_1,..,t_{n}]\otimes \Lambda V,D_2)\ar[u]\ar[r]&
 (\Q [t_1,..,t_{m}]\otimes \Lambda V,D_1)\ar[u]\\
(\Q [t_{m+1},..,t_{n}],0)\ar@{=}[u]\ar[r]&(\Q [t_{1},..,t_{n}],0)\ar[u]\ar[r]&(\Q [t_{1},..,t_{m}],0)\ar[u].
}$$
for $M(BT^{n-m})=(\Q [t_{m+1},..,t_{n}],0)$, 
$M(Y_1)= (\Q [t_1,..,t_{m}]\otimes \Lambda V,D_1)$ and  
$M(Y_2)= (\Q [t_1,..,t_{n}]\otimes \Lambda V,D_2)$.
Then  we simply  write $[D_1]<[D_2]$.  

Even if $r_0(X)>i$
and $\dim H^*(\Q[t_1,\dots,t_{i-1}] \otimes \wedge{V},D)<\infty$,
we may not be able to construct  the KS extension
$$
(\Q[t_i],0)
 \to (\Q[t_1,\dots,t_i] \otimes \wedge{V},D')
 \to (\Q[t_1,\dots,t_{i-1}] \otimes \wedge{V},D)
$$
satisfying $\dim H^*(\Q[t_1,\dots,t_i] \otimes \wedge{V},D')<\infty$
 in general (see Claim 1.2). 

\section{Examples of $r_0(X)\leq 4$}

Refer the arguments of \cite[7.3.2]{FOT} or \cite{JL}
for the computations of toral ranks with minimal models.
We put $M(X)=(\Lambda V,d)$.
A manner to draw $\tilde{\phi}G{\mathcal H}(X)$ often is the following steps.

i) Estimate  $r_0(X)$ by Proposition 2.1.

ii)  Dot $V=\{ (s,t)\in \Z_{\geq 0}\times \Z_{\geq 0}| s\geq 0,t>0,s+t\leq r_0(X) \}$.

iii)  Check whether or not a point `$P$' exists  so that 
$d.i.(P)=(s,r_0(X)-s)\in V$ for $s=1,..,r_0(X)-1$.
If exists (then we say it a {\it bud}), next check the below $\cdots$.  See Claim  1.4 (2).

iv) Check whether or not an edge `$<$' exists  between $P$ and $P'$
with $d.i.(P)=(s,t)$
and $d.i.(P')=(s',t+1)$
for   $s<s'$. See Claim  1.4 (1).
In particular,  the {\it  trunk} 
$P_0-P_1-\cdots  -P_{r_0(X)}$
always exists.

\begin{exmp}{\rm When 
$X=S^{2m+1}\times S^{2n+1}$,
the Hasse diagram of ${\mathcal T}_0(X)$ is totally ordered  as $(1)$
for any $m$ and $n$.
Next put 
 $M(X)=(\Lambda (v_1,v_2,v_3,v_4),d)$
with $dv_1=dv_2=dv_4=0$, $dv_3=v_1v_2$ and 
$|v_1|=|v_2|=3$, $|v_3|=5$, $|v_4|=9$.
It is given by the total space of a non-trival fibration 
$S^5\to X\to S^3\times S^3\times S^9$. 
Then  ${\mathcal H}(X)$ is  given as $(2)$:
$$(1)\ \ \ \ 
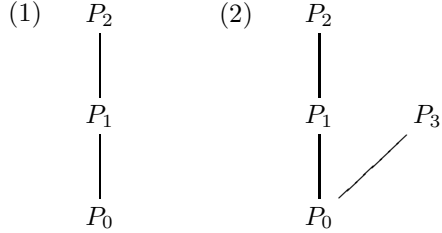
\xymatrix{ 
P_2 \\
P_1 \ar@{-}[u] \\
P_0\ar@{-}[u]
}\ \ \ \ \ \ \ \ \ \ \ 
(2)\ \ \ \ \ 
\xymatrix{ 
P_2 &\\
P_1 \ar@{-}[u] & P_3\\
P_0\ar@{-}[u]\ar@{-}[ur]
}$$
, where $P_1=[ (\Q [t]\otimes \Lambda V,D)]$ 
with $Dv_1=Dv_2=Dv_4=0$ 
$Dv_3=v_1v_2+t^3$, 
 $P_2=[ (\Q [t_1,t_2]\otimes \Lambda V,D)]$ 
with  $Dv_1=Dv_2=0$ 
$Dv_3=v_1v_2+t^3_1$,
$Dv_4=t^5_2$,
  $P_3=[ (\Q [t]\otimes \Lambda V,D)]$ 
with  $Dv_1=Dv_2=0$, 
$Dv_3=v_1v_2$,
$Dv_4=v_1v_3t+t^5$. 
Note that  $\dim H^*(\Q [t,t']\otimes \Lambda V,D')=\infty$
for any KS extension 
$(\Q [t,t']\otimes \Lambda V,D')$ of it. }\end{exmp}

In general,   ${\mathcal H}(X)$ is 
given  as only  $(1)$ or $(2)$ if $r_0(X)=2$.
Next we will consider the cases of $r_0(X)>2$.

A minimal model $(\Lambda V,d)$ is said to be {\it pure}
 if $dV^{even}=0$ and $dV^{odd}\subset \Lambda V^{even}$.
 
\begin{lem}For $m<n$, if 
 $M(Y)=(\Lambda (u_1,..,u_{m},v_1,..,v_n),d)$ 
with $|u_i|$ even and $|v_1|=\cdots =|v_n|$ odd
is a pure model, then
$r_0(Y)=n-m$.
\end{lem}

\noindent
{\it Proof.}
From \cite[Theorem 1]{AH},
 $r_0(Y)\leq n-m$.
From 
\cite[Lemma 8]{H2},
there is a sub-basis $v_1',..,v_m'$ of $\Q(v_1,..,v_n)$
such that $dv_1',..,dv_m'$ is a regular sequence, i.e., 
$\dim H^*(\Lambda (u_1,..,u_{m},v_1',..,v_m'),d)<\infty$.
Then there is a sub-basis
$v_{i_1},..,v_{i_{n-m}}$
with $\Q(v_{i_1},..,v_{i_{n-m}})\oplus \Q(v_1',..,v_m')=\Q(v_1,..,v_n)$.
For $j=1,..,n-m$, put $Dv_{i_j}=dv_{i_j}+t_j^{a_j}$ with $a_j=(|v_{i_j}|+1)/2$.
Then $[t_j^{a_j}]\in \Q [u_1,..,u_m]/(dv_1',..,dv_m')$ and especially
$\dim H^*(\Q [t_{i_1},..,t_{i_{n-m}}]\otimes \Lambda (u_1,..,u_{m},v_1,..,v_n),D)<\infty$.
From Proposition 2.1,  $r_0(Y)\geq n-m$.
\hfill\qed

\begin{thm}
If $X$ has the rational homotopy type of product of odd spheres with  same dimensions,
 $X_{\Q}\simeq (S^k\times \cdots \times S^k)_{\Q}$ for some $k>1$, then 
${\mathcal T}_0(X)$ is totally ordered.
\end{thm}

\noindent
{\it Proof.} Put $r_0(X)=n$.
Suppose $A=(\Q [t_1,..,t_{n-s}]\otimes \Lambda (v_1,..,v_n),D)$
satisfies $\dim H^*(A)<\infty$.
It is easy to check that 
$A$ is pure.
From the above lemma,
$r_0(A)=s$. Thus 
there is no point  $P$ in ${\mathcal H}(X)$ such that
$d.i.(P)=(s,n-s)$ for $s>0$. 
We have done from Claim 1.4 (2).
\hfill\qed

\begin{thm}
Suppose that $1<n_1\leq n_2\leq n_3\leq n_4$ are odd.

(1) For 
$X=S^{n_1}\times S^{n_2}\times S^{n_3}$,
there exists an element  $P$ in ${\mathcal T}_0(X)$
with $d.i.(P)=(2,1)$  
if and only if $n_1+n_2<n_3$.
 
(2) For 
$X=S^{n_1}\times S^{n_2}\times S^{n_3}\times S^{n_4}$,
there exists an element $P$ in ${\mathcal T}_0(X)$
with $d.i.(P)=(3,1)$  if and only if $n_1+n_2<n_3$
and  $n_1+n_3<n_4$.
\end{thm}

\noindent
{\it Proof.} (1) Put $M(X)=(\Lambda (v_1,v_2,v_3),0)$ with $|v_i|=n_i$.
Then 
$Dv_1=Dv_2=0$, $Dv_3=v_1v_2t^{(n_3-n_1-n_2+1)/2}+t^{(n_3+1)/2}$
if and only if $d.i.([D])=(2,1)$.

\noindent
(2) Put $M(X)=(\Lambda (v_1,v_2,v_3,v_4),0)$.
Then there is a differential $D$ with 
$Dv_3=v_1v_2t^{(n_3-n_1-n_2+1)/2}$ and $Dv_4=v_1v_3t^{(n_4-n_1-n_3+1)/2}+t^{(n_4+1)/2}$
if and only if 
there exists a bud $P$ of  $d.i.(P)=(3,1)$.
\hfill\qed

\begin{exmp}{\rm 
Let $X$ be the product of three odd-spheres.
Then $r_0(X)=3$.
From Theorem 3.3, ${\mathcal T}_0(S^3\times S^3\times S^3)$ is  given by
 the 4-points $\{ P_0,P_1,P_2,P_3\}$ 
 which is totally ordered as $(1)$. But
 ${\mathcal T}_0(S^3\times S^3\times S^7)$ 
{\it contains}  the 5-points $\{ P_0,P_1,P_2,P_3,P_4\}$ which is partially ordered as $(2)$:
$$(1) \ \ \ \  
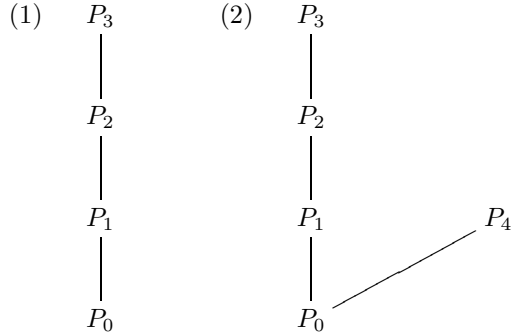
\xymatrix{ 
P_3& \\
P_2 \ar@{-}[u]&\\
P_1 \ar@{-}[u] & \\
P_0\ar@{-}[u]
} \ \ 
(2)\ \ \ \ 
\xymatrix{ 
P_3& &\\
P_2 \ar@{-}[u]&&\\
P_1 \ar@{-}[u] & &P_4\\
P_0\ar@{-}[u]\ar@{-}[urr]
}$$
, where  $P_1=[(S^2\times S^3\times S^7)_{\Q}](=[(S^3\times S^3\times \C P^3)_{\Q}])$,
$P_2=[(S^2\times S^2\times S^7)_{\Q}](=[ (S^2\times S^3\times \C P^3)_{\Q}])$
and $P_3=[(S^2\times S^2\times  \C P^3)_{\Q}]$. 
Here $P_4=[Y_{\Q}]$ is given by 
  $M(Y)=(\Q [t]\otimes \Lambda (x,y,z),D)$ with $Dx=Dy=0$ and $Dz=xyt+t^4$
for $M(X)=(\Lambda (x,y,z),0)$ of $|x|=|y|=3$ and $|z|=7$.
Then $H^*(Y;\Q)\cong \Lambda (x,y)\otimes \Q[t]/(xyt+t^4)$, which is finite dimensional.

Note  $r_0(Y)=0$ from Proposition 2.1.
Indeed, suppose that there is  a KS extension 
$$(\Q [t_2],0)\to (\Q [t_1,t_2]\otimes \Lambda (x,y,z),D')\to (\Q [t_1]\otimes \Lambda (x,y,z),D)=M(Y).$$  We have $D'\circ D'\neq 0$ for any  non-trivial differential 
 $D'x=f(t_1,t_2)$ and  $D'y=g(t_1,t_2)$ in $\Q [t_1,t_2]$.
Also if $D'x=D'y=0$ and $D'z=xyt_1+t_1^4+axyt_2+\sum a_{ij}t_1^{i}t_2^j$ ($a, a_{ij}\in \Q$),
then $\dim H^*(\Q [t_1,t_2]\otimes \Lambda (x,y,z),D')=\infty$
for any $a,\ a_{ij}$.
Thus $d.i.(P_4)=(2,2-0-1)=(2,1)$. 

Thus the set  ${\mathcal T}_0(X)$
is more sensitive than the number 
$r_0(X)$ about degrees of the  rational homotopy group of $X$.
}
\end{exmp}


\begin{exmp}{\rm 
Put $M(X)=(\Lambda (v_1,v_2,v_3,v_4,v_5),d)$
with $dv_1=dv_2=dv_4=dv_5=0$, $dv_3=v_1v_2$ and 
$|v_1|=|v_2|=3$, $|v_3|=5$, $|v_4|=9$, $|v_5|=15$. Then
 ${\mathcal H}(X)$ is  given as
$$
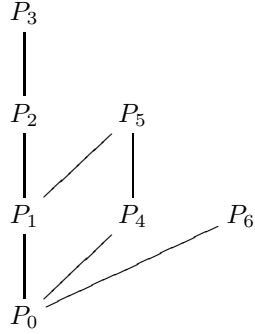
\xymatrix{ 
P_3 &  &   \\
P_2 \ar@{-}[u]& P_5&  \\
P_1 \ar@{-}[u]\ar@{-}[ru] & P_4\ar@{-}[u]& P_6 \\
P_0\ar@{-}[u]\ar@{-}[ur]\ar@{-}[urr]\\
}\
$$
, where $P_1=[ (\Q [t]\otimes \Lambda V,D)]$ 
with $Dv_1=Dv_2=Dv_5=0$, $Dv_4=t^5$,  
$Dv_3=v_1v_2$.
  $P_2=[ (\Q [t_1,t_2]\otimes \Lambda V,D)]$ 
with  $Dv_1=Dv_2=Dv_5=0$, 
$Dv_3=v_1v_2+t^3_2$,
$Dv_4=t^5_1$.\\
  $P_3=[ (\Q [t_1,t_2,t_3]\otimes \Lambda V,D)]$ 
with  $Dv_1=Dv_2=0$, 
$Dv_3=v_1v_2+t_2^3$,
$Dv_4=t_1^5$, $Dv_5=t_3^8$.\\
$P_4=[ (\Q [t]\otimes \Lambda V,D)]$ 
with  $Dv_1=Dv_2=Dv_5=0$,  $Dv_3=v_1v_2$,
$Dv_4=v_1v_3t+t^5$.\\
 $P_5=[ (\Q [t_1,t_2]\otimes \Lambda V,D)]$ 
with  $Dv_1=Dv_2=0$,  $Dv_3=v_1v_2$, 
$Dv_4=v_1v_3t_1+t^5_1$, 
$Dv_5=t^8_2$.\\
$P_6=[ (\Q [t]\otimes \Lambda V,D)]$ 
with  $Dv_1=Dv_2=Dv_4=0$,   $Dv_3=v_1v_2$, 
$Dv_5=v_1v_3t^4+v_2v_4t^2+t^8$.
}
\end{exmp}

\begin{rem}{
If $r_0(X)=3$,  ${\mathcal H}(X)$ is  given as
(1), (2) of Example 3.5, Example 3.6 or as one of the five diagrams:
$$
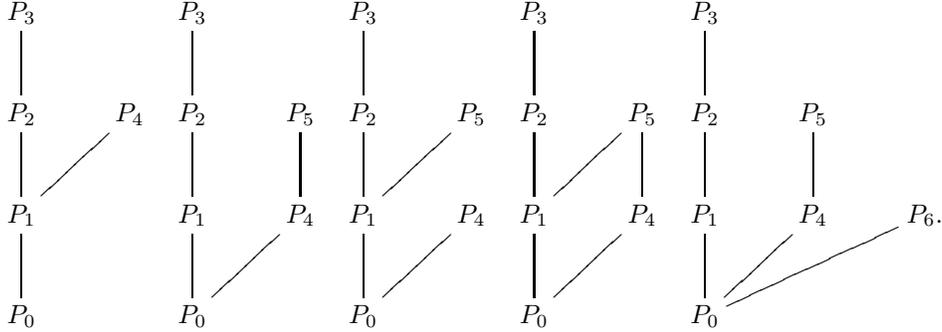
\xymatrix{ 
P_3& \\
P_2 \ar@{-}[u]& P_4\\
P_1 \ar@{-}[u]  \ar@{-}[ur]& \\
P_0\ar@{-}[u]
}
\ \ 
\xymatrix{ 
P_3& \\
P_2 \ar@{-}[u]& P_5\\
P_1 \ar@{-}[u] & P_4 \ar@{-}[u]\\
P_0\ar@{-}[u]\ar@{-}[ur]
}
\ \ 
\xymatrix{ 
P_3& \\
P_2 \ar@{-}[u]& P_5\\
P_1 \ar@{-}[u]  \ar@{-}[ur]& P_4\\
P_0\ar@{-}[u]\ar@{-}[ur]
}
\ \ 
\xymatrix{ 
P_3& \\
P_2 \ar@{-}[u]& P_5\\
P_1 \ar@{-}[u]  \ar@{-}[ur] & P_4 \ar@{-}[u]\\
P_0\ar@{-}[u]\ar@{-}[ur]
}
\ \ 
\xymatrix{ 
P_3& &\\
P_2 \ar@{-}[u]& P_5 &\\
P_1 \ar@{-}[u]   & P_4 \ar@{-}[u] & P_6.\\
P_0\ar@{-}[u]\ar@{-}[ur]\ar@{-}[urr]
}
$$
Thus we see that $2\leq \sharp \{ {\mathcal H}(X)|r_0(X)=3\}\leq 8$.
}
\end{rem}

Finally we give two examples with the same Hasse diagrams.

\begin{exmp}{\rm 
Put $M(X_1)=(\Lambda (v_1,v_2,v_3,v_4,v_5,v_6),d)$
with $dv_1=dv_2=dv_4=dv_5=dv_6=0$, $dv_3=v_1v_2$ and 
$|v_1|=|v_2|=3$, $|v_3|=5$, $|v_4|=9$, $|v_5|=13$, $|v_6|=17$. Then
 ${\mathcal H}(X_1)$ is  given as 

$$
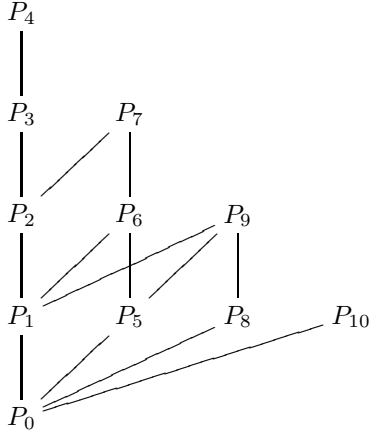
\xymatrix{ 
P_4 & &  & \\
P_3 \ar@{-}[u]& P_7 &  & \\
P_2 \ar@{-}[u]\ar@{-}[ru]& P_6 \ar@{-}[u]& P_{9}& \\
P_1 \ar@{-}[u]\ar@{-}[ru]\ar@{-}[urr] & P_5\ar@{-}[u]\ar@{-}[ru] & P_8\ar@{-}[u]& P_{10} \\
P_0\ar@{-}[u]\ar@{-}[ur]\ar@{-}[urr]\ar@{-}[urrr]\\
}$$
, where
 $P_4=[(\Q[t] \otimes \wedge{V},D)]$ with $Dv_1=Dv_2=0$,
$Dv_3=v_1v_2+t^3$, $Dv_4=t^5$, $Dv_5=t^7$, $Dv_6=t^9$.\\
 $P_5=[(\Q[t] \otimes \wedge{V},D)]$ with $Dv_1=Dv_2=Dv_5=Dv_6=0$,
$Dv_3=v_1v_2$, $Dv_4=v_1v_3t+t^5$.\\
$P_8=[(\Q[t] \otimes \wedge{V},D)]$ with 
$Dv_1=Dv_2=Dv_4=Dv_6=0$, $Dv_3=v_1v_2$, 
$Dv_5=v_2v_4t+v_1v_3t^4+t^7$.\\
$P_9=[(\Q[t_1,t_2] \otimes \wedge{V},D)]=[(\Q[t_1,t_2] \otimes \wedge{V},D')]=[(\Q[t_1,t_2] \otimes \wedge{V},D'')]$ with 
$Dv_1=Dv_2=0$, $Dv_3=v_1v_2$ (same for $D'$ and $D''$) and \\
$$Dv_4=v_1v_3t_2+t_1^5,  \ Dv_5=0, \ Dv_6=v_2v_5t_2+t_2^9,$$
$$D'v_4=v_1v_3t_1+t_1^5,  \ D'v_5=0,    \ D'v_6=v_2v_5t_2+t^9_2,$$
$$D''v_4=v_1v_3t_1+t_2^5,  \ D''v_5=0, \ D''v_6=v_2v_5t_1+t_1^9.$$
Note that $P_1<P_9$ is given by $D$, $P_5<P_9$ by $D'$
and  $P_8<P_9$  by $D''$.\\
$P_{10}=[(\Q[t] \otimes \wedge{V},D)]$ with $Dv_1=Dv_2=Dv_4=Dv_5=0$, $Dv_3=v_1v_2$
and
$Dv_6=v_1v_5t+v_2v_4t^3+v_1v_3t^5+t^{9}$.

}
\end{exmp}

\begin{exmp} {\rm
Put $M(X_2)=(\Lambda (v_1,v_2,v_3,v_4,v_5,v_6,w_1,w_2,w_3,w_4,w_5,w_6,
u,y,z,v),d)$ with
$dv_1=\cdots =dv_6=dw_1=\cdots =dw_6=du=dy=dz=0$,  $$dv=v_1v_2v_3v_4+v_1v_2v_5v_6+w_1w_2w_3w_4+w_1w_2w_5w_6+u^2$$
and $|v_1|=\cdots =|v_6|=|w_1|=\cdots =|w_6|=3$, $|u|=6$,  $|y|=|z|=7$, $|v|=11$.
We see  $\dim H^*(X_2;\Q )<\infty$
since 
$M(X_2)$
is the total space of a KS extension
$$(\Lambda (v_1,\cdots ,v_6,w_1,\cdots ,w_6,y,z),0)\to 
(\Lambda V,d)\to (\Lambda (u,v),\overline{d})=M(S^6),$$
where $\overline{d}u=0$ and $\overline{d}v=u^2$.
Remark  the space $X_1$
in Example 3.8  is not a formal space  but  $X_2$ is  formal  \cite{FHT}.
Note that  $r_0(X')=0$ for a space $X'$ with ${X_2}_{\Q}\simeq (X'\times S^7\times S^7)_{\Q}$, 
where $M(X')=(\Lambda (v_1,v_2,v_3,v_4,v_5,v_6,w_1,w_2,w_3,w_4,w_5,w_6,
u,v),d)$.
  
Then  ${\mathcal H}(X_2)$ is   given as 

$$
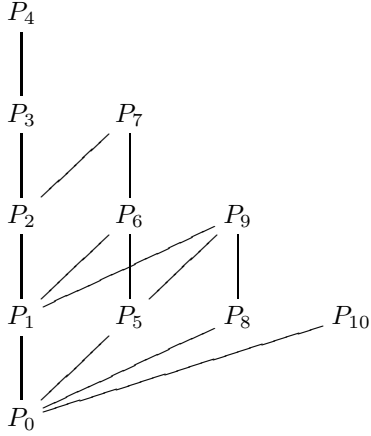
\xymatrix{ 
P_4 & &  & \\
P_3 \ar@{-}[u]& P_7 &  & \\
P_2 \ar@{-}[u]\ar@{-}[ru]& P_6 \ar@{-}[u]& P_{9}& \\
P_1 \ar@{-}[u]\ar@{-}[ru]\ar@{-}[urr] & P_5\ar@{-}[u]\ar@{-}[ru] & P_8\ar@{-}[u]& P_{10} \\
P_0\ar@{-}[u]\ar@{-}[ur]\ar@{-}[urr]\ar@{-}[urrr]\\
}$$
, where  $P_1=[ (\Q [t]\otimes \Lambda V,D)]$ 
such that $Dv_i=0$ for  $i\neq 4$, 
$Dv_4=t^2$,  $Dw_i=Du=Dz=0$,
 $Dy=v_2v_3t$,
 $Dv=dv-v_1yt$. \\
$P_4=[  (\Q [t_1,t_2,t_3,t_4]\otimes \Lambda V,D)]$ with  
$$Dv_1=Dv_2=Dv_3=Dv_5=Dw_1=Dw_2=Dw_3=Dw_5=Du=0,\cdots (*)$$ 
$$Dv_4=t_1^2,\ Dv_6=t_2^2, \ Dw_4=t_3^2, \ Dw_6=t_4^2,$$ 
 $$Dy=v_2v_3t_1+v_1v_5t_2,
\ Dz=w_2w_3t_3+w_1w_5t_4,$$
 $$Dv=dv-v_1yt_1+v_2yt_2-w_1zt_3+w_2zt_4.$$
 Then  $D\circ D=0$
and $\dim H^*(\Q [t_1,t_2,t_3,t_4]\otimes \Lambda V,D)<\infty$.
Thus $r_0(X_2)\geq 4$ and we deduce  $r_0(X_2)<5$
by the direct (but complicated) calculations that
$\dim H^*(\Q [t_1,t_2,t_3,t_4,t_5]\otimes \Lambda V,D)=\infty$
for any $D$. 
Note the part $(*)$ of $P_4$  is applied  for all differentials below.\\
$P_5=[ (\Q [t]\otimes \Lambda V,D)]$ 
with $Dv=dv$, $Dy=t^4$, $Dv_i=Dw_i=Du=Dz=0$.
\\
$P_6=[ (\Q [t_1,t_2]\otimes \Lambda V,D)]$ 
such that $Dv_i=0$ for  $i\neq 4$, 
$Dv_4=t_1^2$,  $Dw_i=Du=0$,
 $Dy=v_2v_3t_1$,
 $Dv=dv-v_1yt_1$, $Dz=t_2^4$. \\
$P_8= [ (\Q [t]\otimes \Lambda V,D)]$
with $Dv=dv-v_1yt+v_2yt-w_1zt$, $Dy=v_2v_3t+v_1v_5t$,
$Dz=w_2w_3t$, $Dv_4=Dv_6=Dw_4=t^2$, $Dw_6=0$.\\
 $P_9=[ (\Q [t_1,t_2]\otimes \Lambda V,D)]=[ (\Q [t_1,t_2]\otimes \Lambda V,D')]=[ (\Q [t_1,t_2]\otimes \Lambda V,D'')]$ with
$$Dv_4=t_1^2,\ Dv_6=Dw_4=Dw_6=t_2^2,\ Dy=v_2v_3t_1+v_1v_5t_2,$$
$$ Dz=w_2w_3t_2+w_1w_5t_2,\ 
Dv=dv-v_1yt_1+v_2yt_2-w_1zt_2+w_2zt_2,$$
 $$D'v_i=D'w_i=D'u=0, \ D'y=t_1^4, \ D'z=t_2^4,\ D'v=dv,$$
$$D''v_4=t_2^2,\ D''v_6=D''w_4=D''w_6=t_1^2,\ D''y=v_2v_3t_2+v_1v_5t_1,$$
$$ D''z=w_2w_3t_1+w_1w_5t_1,\ 
D''v=dv-v_1yt_2+v_2yt_1-w_1zt_1+w_2zt_1.$$
Note that $P_1<P_9$ is given by $D$, $P_5<P_9$  by $D'$
and  $P_8<P_9$  by $D''$.\\
 $P_{10}= [ (\Q [t]\otimes \Lambda V,D)]$
with $Dv_4=Dv_6=Dw_4=Dw_6=t^2$, $Dy=v_2v_3t+v_1v_5t$,
$Dz=w_2w_3t+w_1w_5t$,  $Dv=dv-v_1yt+v_2yt-w_1zt+w_2zt$.
}
\end{exmp}

\section{Proof of Theorem 1.5}

Let $G$ be a connected, non-directed,  finite, simple(i.e., without multiple edges, 
loops),  based graph with the vertex set $V(G)=\{v_0,v_1,..,v_N\}$ of the base point $v_0$.
For 
the set of distances $D_0=\{d(v_0,v_i)|v_i\in V(G)\}_i$ 
between the points of $G$ and  $v_0$,
put $n=max D_0$. 
Suppose that  a path  of  length  $n$
$$ l_0\ :\ \ \ v_{0}\to v_{i_1}\to \cdots \to v_{i_{n-1}}\to v_{i_n} $$
with  $d(v_0,v_{i_n})=n$ is unique.\ \ \ \ \  $\cdots (0)$\\
Then put  $\psi (v_0):=(0,0)$ and 
$$\psi (v_{i_u}):=(0,u)\in \Z_{\geq 0}\times \Z_{\geq 0}$$
for $u=1,..,n$.

Next 
put $D_1=\{d(v_0,v_j)|v_j\in V_1=V(G)-V(l_0) \}_j$ and  
the set of paths with  length   $n_1=max D_1$ as
$$L_1=\{ l_{1,j}\}_j=\{  v_{0}\to v_{j_1}\to \cdots \to v_{j_{n_1-1}}\to v_{j_{n_1}}|\ d(v_0,v_{j_{n_1}})=n_1, v_{j_{n_1}}\in V_1\}_j.$$
Here $V(l_0) =\{ v_0,v_{i_1},.., v_{i_n} \} $.
Suppose that (for some $c$)
$$j_m\neq i_m\mbox{ for }m>c\mbox{
 if }j_c\neq i_c.\ \ \  \cdots (1)$$
For a  path $l_{1,j}$ of  $L_1$, if $v_{j_c}=v_{i_c}$ for $c=0,..,m-1$ and $v_{j_m}\neq v_{i_m}$,
put $$\psi (v_{j_u}):=(n-n_1,u)\in  \Z_{\geq 0}\times \Z_{\geq 0}$$
for $u=m,m+1,..,n_1$.

Next put $D_2=\{d(v_0,v_k)|v_k\in V_2=V(G)-(V(l_0) \cup V(L_1)) \}_k$
and  
the set of paths of  with length $n_2=max D_2$ 
$$L_2=\{ l_{2,k}\}_k=\{  v_{0}\to v_{k_1}\to \cdots \to v_{k_{n_2-1}}\to v_{k_{n_2}}|\ d(v_0,v_{k_{n_2}})=n_2, v_{k_{n_2}}\in V_2\}_k.$$
Suppose that (for some $c$)
$$k_m\neq i_m \mbox{ for }m>c
\mbox{ if }k_c\neq i_c\mbox{ \ \ and}$$  
$$k_m\neq i_m,j_m \mbox{ for }m>c 
\mbox{ if }k_c\neq j_c.\ \ \   \cdots (2)$$
For a path $l_{2,k}$ of  $L_2$, if $v_{k_c}=v_{i_c}$ or $v_{k_c}=v_{j_c}$
 for $c=0,..,m-1$ but   $v_{k_m}\neq v_{i_m}$ and $v_{k_m}\neq v_{j_m}$,
put $$\psi (v_{k_u}):=(n-n_2,u)\in  \Z_{\geq 0}\times \Z_{\geq 0}$$
for $u=m,m+1,..,n_2$.

Iterating this argument, we have an injection $\psi :V(G)\to  \Z_{\geq 0}\times \Z_{\geq 0}$
and it is naturally extended to the map from the set  of 
edges, $\tilde{\psi} :E(G)\to  \R_{\geq 0}\times \R_{\geq 0}$  as  $\tilde{\psi}(v_av_b)=\psi (v_a)-\psi (v_b)$, the line segment 
with extremal points 
$\psi (v_a)$ and $\psi (v_b)$, for any edge $v_av_b$ of $G$.
Thus there is 
the embedding of $G$ into $  \R_{\geq 0}\times \R_{\geq 0}$

$$\xymatrix{ 
V(G)\ar[r]^{\psi\  \ \ }\ar[d]_{\cap}& \Z_{\geq 0}\times \Z_{\geq 0}\ar[d]^{\cap}\\
G\ar[r]^{\tilde{\psi}\ \ \ }&  \R_{\geq 0}\times \R_{\geq 0}.\\
}
$$

Notice that two graphs 
$G$ and $G'$ satisfying (0),(1),(2),.. are isomorphic as based graphs
if and only if $\tilde{\psi}G= \tilde{\psi}G'$.\\

\noindent
{\it Proof of Theorem 1.5}.
 If $G=G{\mathcal H}(X)$, the above conditions (0), (1), (2),.. are satisfied
from Lemma 1.3.
 Thus the above map $\tilde{\psi}$ is defined and  we see
 $\tilde{\psi} G{\mathcal H}(X)=\tilde{\phi} G{\mathcal H}(X)$
in $\R_{\geq 0}\times \R_{\geq 0}$. 
Suppose  ${\mathcal H}(X)\neq {\mathcal H}(Y)$.
 Then 
$\tilde{\psi} G{\mathcal H}(X)=\tilde{\phi} G{\mathcal H}(X)\neq \tilde{\phi} G{\mathcal H}(Y)=
\tilde{\psi} G{\mathcal H}(Y)$.
Thus
 $G{\mathcal H}(X)$ and $G{\mathcal H}(Y)$ are not isomorphic  as based graphs. 
\hfill\qed

\section{Appendix}

Recall an edge  of $G{\mathcal H}(X)$
is represented  by a rationalized Borel fibration 
$$Y_{\Q}\to (ES^1\times_{S^1}Y)_{\Q}\to BS^1_{\Q}$$
where $Y_{\Q} \in {\mathcal X}_n$
and $(ES^1\times_{S^1}Y)_{\Q} \in {\mathcal X}_{n+1}$ for some $n$.
It is given as 
$$\xymatrix{ 
\bullet\\
\bullet \ar@{-}[u]
}\ \ \ \ {or} \ 
\xymatrix{ 
&\bullet\\
\bullet \ar@{-}[ru]&
}
\ \ \  {or} \ \xymatrix{ 
&&\bullet\\
\bullet \ar@{-}[rru]&&
} \ {or} \ \ \cdots 
$$
in 
 $\tilde{\phi}G{\mathcal H}(X)$.

\begin{defi}
Suppose $Y_{\Q} \in {\mathcal X}_n$.
For two elements 
$Y_1=(ES^1\times_{S^1}^{\mu_1}Z_1)_{\Q}$ and 
$Y_2=(ES^1\times_{S^1}^{\mu_2}Z_2)_{\Q}$ of ${\mathcal X}_{n+1}$,
we denote 
$$Y_1\underset{Y_3}{\sim}Y_2$$ 
if there exists 
 a  homotopy commutative diagram 
of fiber inclusions of rationalized  Borel fibrations over 
$BS^1_{\Q}$
 $$\xymatrix{  Y_{\Q}\ar[d]\ar[r] & (ES^1_2\times_{S^1_2}Z_2)_{\Q}=Y_2 \ar[d]\\
Y_1=(ES^1_1\times_{S^1_1}Z_1)_{\Q}\ar[r] &(E(S^1_1\times S^1_2)\times_{S^1_1\times S^1_2}Z_3)_{\Q} =:Y_3}$$
 where ${Z_1}_{\Q}\simeq {Z_2}_{\Q}\simeq {Z_3}_{\Q}\simeq Y_{\Q}$
and $\dim H^*(Y_3;\Q )<\infty$. 
\end{defi}

 Note $Y_3\in {\mathcal X}_{n+2}$
and  in general $r_0(Y_1)\neq r_0(Y_2)$.
 The Sullivan model is given as the DGA-homotopy commutative diagram of
natural projections 
 $$\xymatrix{ M(Y)\cong (\Lambda W,d_W)\ \  &
(\Q [t_{2}]\otimes \Lambda W,D_2)\ar[l]\\ 
(\Q [t_1]\otimes \Lambda W,D_1)\ar[u]&(\Q [t_1,t_2]\otimes \Lambda W,D) \ar[u]\ar[l]& }$$
 with $\dim H^*(\Q [t_1,t_2]\otimes \Lambda W,D)<\infty$.
 Here $W=\Q (t_3,..,t_{n+2})\oplus V$ for $M(X)=(\Lambda V,d)$ and $d_W|V=d$.
Remark that Definition 5.1  is not an equivalence relation.

\begin{defi}
For  edges (1-cells) $P_aP_b$, $P_aP_d$, $P_bP_c$ and $P_dP_c$
in $G{\mathcal H}(X)$,
which is given as (a horizontal deformation of) 
$$
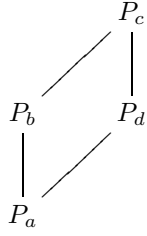
\xymatrix{ 
&  P_c \\
P_b \ar@{-}[ur] & P_d\ar@{-}[u]\\
P_a\ar@{-}[u]\ar@{-}[ur]&\\
}$$
in 
 $\tilde{\phi}G{\mathcal H}(X)$, 
we say that a 2-cell attachs on the 1-cycle  $ \square{P_aP_bP_cP_d}$ 
(or  simply that  $ \square{P_aP_bP_cP_d}$ makes a leaf)
and denote as
$\partial e^2=\square{P_aP_bP_cP_d}$ 
if  $Y_1\sim_{Y_3}Y_2$ for $[Y_{\Q}]=P_a$, $[Y_1]=P_b$, $[Y_2]=P_d$
and $[Y_3]=P_c$.
\end{defi}

The existence of a leaf may  
depend  on the  degree of certain freedom of $\{ D\}$
 that represent  the upper right point $P_c$ of a cycle $ \square{P_aP_bP_cP_d}$.
In Example 3.8,
we easily find that  
$ \square{P_0P_5P_9P_8}$
makes a leaf
by $Dv_1=Dv_2=Dv_4=0$, $Dv_3=v_1v_2$,
$$Dv_5=v_1v_4t_1+t_1^7\ \ \mbox{ and}$$
$$Dv_6=v_2v_4t_2^3+v_1v_3t_2^5+t_2^9,$$
where $[D]=P_9$, $[D_1]=P_5$ and $[D_2]=P_8$.
Indeed,  then the above DGA-diagram is commutative. 
But,  in Example 3.9, 
the author can not find a differential  $D$ that makes 
the above homtopy commutative diagram for 
the 1-cycle $ \square{P_0P_5P_9P_8}$.

In general, if $G{\mathcal H}(X)$ contains 
(a horizontal deformation of) 
$$
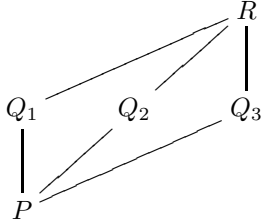
\xymatrix{ 
& & R& \\
Q_1 \ar@{-}[urr] & Q_2\ar@{-}[ru] & Q_3\ar@{-}[u]& \\
P\ar@{-}[u]\ar@{-}[ur]\ar@{-}[urr]\\
}$$
as a sub-graph
with
$\partial e^2_1=\square{PQ_1RQ_2}$,
$\partial e^2_2=\square{PQ_2RQ_3}$
and 
$\partial e^2_3=\square{PQ_1RQ_3}$,
then ${\mathcal K}(X)$ contains 
$$(\ e_1^2\cup e_2^2\ )\cup_{\square{PQ_1RQ_3}} e_3^2\ \ \cong \ \ S^2.$$
Thus three pieces of   leaf
can make a 2-sphere.  

\begin{rem}
To append  certain  further informations of ${\mathcal X}$
on ${\mathcal T}_0(X)$,
it may be suitable  to regard 
(${\mathcal T}_0(X)$
as the 0-skeleton and) the based graph $G{\mathcal H}(X)$  as
the 1-skeleton of a   finite  CW complex ${\mathcal T}(X)$,
which is obtained by generalizing Definition 5.2. 
When  $\tilde{\phi}G{\mathcal H}(X)$
 contains
(a horizontal deformation of) 

$$
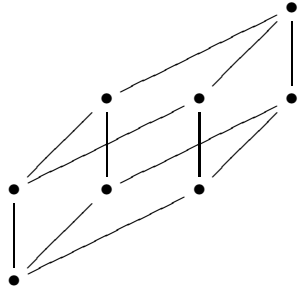
\xymatrix{ 
&&&\bullet\\
&\bullet \ar@{-}[urr]&\bullet \ar@{-}[ur]&\bullet \ar@{-}[u]\\
\bullet \ar@{-}[ur] \ar@{-}[urr]&\bullet \ar@{-}[u] \ar@{-}[urr]&\bullet \ar@{-}[u] \ar@{-}[ur]&\\
\bullet\ar@{-}[u] \ar@{-}[ru] \ar@{-}[rru]&&&
}$$

\noindent
as a sub-graph,
then   
 a 3-cell  of ${\mathcal T}(X)$ is given by  the existence of 
the homotopy commutative digram of natural projections
 $$
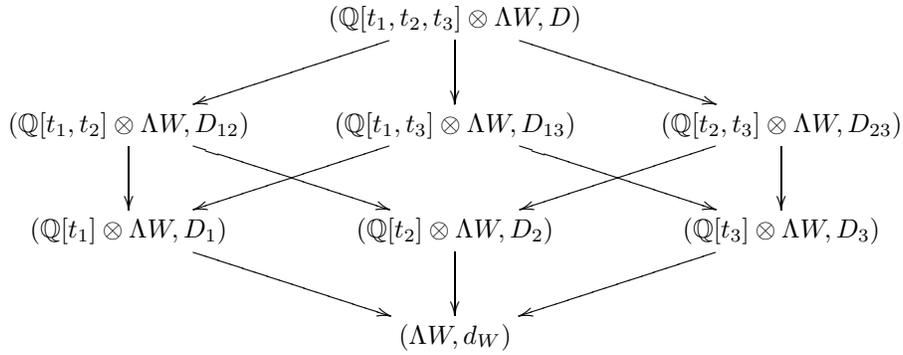
\xymatrix{ &(\Q [t_1,t_2,t_3]\otimes \Lambda W,D)\ar[ld]\ar[d]\ar[dr]&\\
 (\Q [t_1,t_2]\otimes \Lambda W,D_{12})\ar[d]\ar[rd]&(\Q [t_1,t_3]\otimes \Lambda W,D_{13})\ar[ld]\ar[rd]&(\Q [t_2,t_3]\otimes \Lambda W,D_{23})\ar[d]\ar[dl]\\
 (\Q [t_1]\otimes \Lambda W,D_1)\ar[rd]&(\Q [t_2]\otimes \Lambda W,D_2)\ar[d]&(\Q [t_3]\otimes \Lambda W,D_3)\ar[ld]\\
& (\Lambda W,d_W)&
 }$$
which represents the above sub-graph.
Similarly we can construct  higher dimensional CW-structure,
 which makes a complex
 ${\mathcal T}(X)$.
It must be  a topological homotopy invariant of spaces.
 (The complex 
 ${\mathcal T}(X)$ 
is at most 2-dimensional if $r_0(X)\leq 5$.)   
If ${\mathcal T}(X)$ is compared
 to a plant, then
the base point $X_{\Q}$ corresponds to the seed  (that grows up to be the plant),
and $BS^1_{\Q}$, the water (that is necessary for its growth).
\end{rem}

\end{document}